\documentclass[leqno]{article}
\usepackage{epsfig}
\input epsf
\usepackage{latexsym}
\usepackage{amsmath}
\title{ Stable closed equilibria for anisotropic surface energies: Surfaces with edges}
\author{ By  B{\footnotesize ENNETT} P{\footnotesize ALMER}}

\begin{document}

\newtheorem{theorem}{Theorem}[section]
\newtheorem{cor}{Corollary}[section]
\newtheorem{prop}{Proposition}[section]
\newtheorem{lemma}{Lemma}[section]
\newtheorem{condition}{Condition}[section]
\newtheorem{example}{Example}[section]

\newtheorem{definition}{Definition}[section]
\newtheorem{remark}{Remark}[section]
\newtheorem{conjecture}{Conjecture}[section]
\newtheorem{claim}{Claim}[section]
\newtheorem{question}{Question}[section]
\newcommand{\rf}[1]{\mbox{(\ref{#1})}}
\maketitle
\begin{abstract} We study the stability of closed, not necessarily smooth, equilibrium surfaces of an anisotropic surface energy
for which the Wulff shape is not necessarily smooth. We show that if the Cahn Hoffman field
can be extended continuously to the whole surface and if the surface is stable, then the surface is,
up to rescaling, the Wulff shape.
\end{abstract}

In this paper, we will study the stability of closed surfaces which are in equilibrium for an anisotropic surface energy. Neither the equilibrium surface $\Sigma$ nor the Wulff shape $W$ are assumed to be smooth; they may have edges as depicted in Figure 1. The equilibrium conditions for $\Sigma$ are expressed by the anisotropic mean curvature being constant on each face of $\Sigma$ and at each edge, the force balancing condition found by Cahn and Hoffman \cite{CH} is satisfied. The main result states that if the the Cahn-Hoffman field extends to a continuous map of $\Sigma$ to $W$ and if the surface is stable, then the surface agrees with $W$ up to rescaling. Stability means that the second variation of energy is non negative for those variations preserving the volume enclosed by $\Sigma$. The condition of extendability of the Cahn-Hoffman map seems natural to us since it implies the force balancing condition along the edges and we can produce many
(non closed) examples for which it holds.

The first theorem of the type described above is due to Barbosa-do Carmo, \cite{BdC} who showed that the spheres are the unique stable closed constant mean curvature hypersurfaces in Euclidean space. Later, Wente \cite{W} greatly simplified the proof of this result. Wente's idea was used by the author \cite{P}  to extend the results of
\cite{BdC} to the anisotropic case, assuming both the smoothness of the hypersurface and the smoothness of the Wulff shape. A proof of this result more in the style of \cite{BdC} was given in 
\cite{Wi}.  Similar results involving  higher order anisotropic mean curvatures can be  found in \cite{HeLi2}.

Because of the applications of anisotropic surface energies to studying interfaces of structured materials, we have decided to limit our attention to the case of two dimensional surfaces in three space. It is very likely that these ideas can be extended to higher dimensions. In addition, we have limited our discussion for the most part to embedded surfaces which is not essential.

Besides  the Barbosa-do Carmo theorem, several other results characterizing the sphere among constant mean curvature surfaces have recently been generalized to the anisotropic case.
In \cite{ghlm}, it is shown that the unique embedded closed constant anisotropic mean curvature surfaces are the rescalings of the Wulff shape. Also, in \cite{KP11},  it is shown that the unique closed genus zero constant anisotropic mean curvature surfaces in ${\bf R}^3$ are rescalings of the Wulff shape. These two results are , respectively, generalizations of the Alexandrov and Hopf Theorems.  In both cases, the hypersurfaces and the Wulff shape are assumed to be smooth. It would be interesting to know if either of these results could be generalized further to the type of functionals considered in this paper.
\\[4mm]

\begin{figure}
  \hfill
  \begin{minipage}[n]{.45\textwidth}
    \begin{center}
       \framebox{\includegraphics[width=60mm,height=60mm,angle=0]{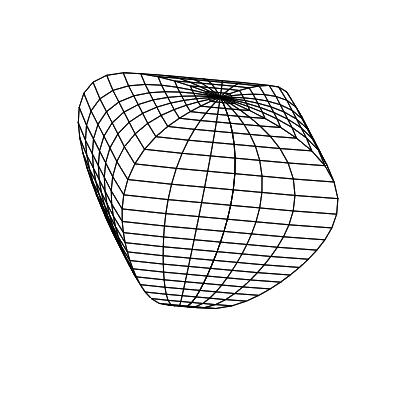}}
    \end{center}
 \end{minipage}
  \hfill
  \begin{minipage}[z]{.45\textwidth}
    \begin{center}
      \framebox{\includegraphics[width=60mm,height=60mm,angle=0]{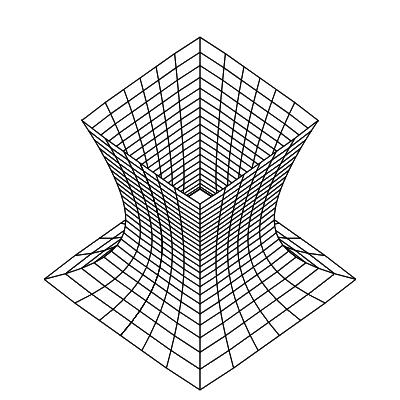}}
    \end{center}
   \end{minipage}
\caption{ A Wulff shape and its equilibrium catenoid}
    \end{figure}

A surface $S\subset {\bf R}^3$ will be called piecewise smooth if there exists a collection $\Gamma_1,...,\Gamma_L$  of pairwise disjoint embedded $C^1$ curves in $S$ such that $S^*:=S\setminus
\Gamma_1\cup...\cup \Gamma_L$ is a disjoint union of smooth open surfaces $S_1,...,S_N$.  In addition, each $S_j$ together with its boundary curves is assumed to be a smooth surface with piecewise smooth 
boundary. The $S_j$'s will be referred to as {\it faces} of $S$ , while the $\Gamma_i$'s will be called edges. 

  We let $\gamma$ be a positive $C^3$ function on $S^2$ and define its Wulff shape $W$ by
  \begin{equation}
\label{W} W=\partial  \bigcap_{n\in S^2} \{ Y\cdot n\le \gamma(n)\}\:,\end{equation}
Despite the smoothness of $\gamma$, $W$ need not be smooth. In this paper, we will assume that $W$ is piecewise smooth in the sense described above.  In addition we will assume that $W$ is convex in the sense that each face of $W$ has uniformly positive curvature $K_W$.

The function $\gamma$ defines an anisotropic surface energy density. If $\Sigma$ is an oriented piecewise smooth surface, then the normal field $\nu$ of $\Sigma$ is almost everywhere defined and the energy of $\Sigma$ is
$${\cal F}[\Sigma]=\int_\Sigma \gamma(\nu)\:d\Sigma\:.$$
(Throughout this paper, we adopt the usual conventions of Lebesgue integration so that a function which is only defined almost everywhere but is continuous and bounded on its domain
is integrable.)

The Cahn-Hoffman field  is defined as follows. Let ${\tilde \gamma}$ denote the positive, degree one homogeneous extension of $\gamma$ and let $\nabla {\tilde \gamma}$ denote its gradient on ${\bf R}^3\setminus \{0\}.$
\begin{equation}
\label{CH}\xi(p) := \nabla {\tilde \gamma}(\nu(p))=D\gamma_{\nu(p)}+\gamma(p)\nu\:, \:p\in \Sigma ^*\:,\end{equation} 
where $D\gamma$ denotes the gradient of $\gamma$ on $S^2$. 

The first variation formula gives
\begin{equation}
\label{deltaf}\delta {\cal F}=-\sum_j\bigl(\int_{\Sigma_j} \Lambda \delta X\cdot \nu \:d\Sigma -\oint_{\partial \Sigma_j} (\xi \times  dX)\cdot \delta X\bigr)\:.
\end{equation}
If $V$ denotes the (algebraic) volume enclosed by the surface, then its variation is
$$\delta V= 3\int_{\Sigma^*} \delta X\cdot \nu\:d\Sigma\:.$$

 \begin{prop}
 A piecewise smooth surface $\Sigma$ is in equilibrium for the volume constrained energy functional if and only if there holds:
 (i) $\Lambda \equiv $ constant on $\cup \Sigma_j\:,$ and (ii)
 $ (\xi_1-\xi_2)\times \Gamma '\equiv 0$ on $\Gamma_i$.
 
 In particular, (ii) holds if $\xi$ extends continuously to all of $\Sigma$. 
 \end{prop}

For a small displacement on a surface represented by a vector ${\vec l}$, Cahn and Hoffman interpreted $\xi \times {\vec l}$ as the force of anisotropic surface tension acting on ${\vec l}$.
If two faces ${\Sigma_1}$,${\Sigma_2} $ meet along the edge  ${\Gamma}$, then with respect to the orientations of these faces,
$\Gamma$ must be counted with one orientation for $\Sigma_1$ and with the opposite orientation for $\Sigma_2$.  Thus $(\xi_1-\xi_2)\times \Gamma' =0$ means that the forces balance along the direction of the edge.

We will call an equilibrium surface {\it stable} if for any volume preserving deformation
$X_\epsilon$ of $X$ through piecewise smooth immersion which is at least twice differentiable with respect to $\epsilon$, we have
$\partial^2_{\epsilon \epsilon}{\cal F}[X_\epsilon]_{\epsilon=0} \ge 0 $ .
\begin{theorem}
\label{th1} Let $\Sigma$ be a closed embedded piecewise smooth surface which is in equilibrium for an anisotropic surface energy having piecewise smooth convex Wulff shape $W$ and assume that the Cahn-Hofmann map extends to a continuous map
$\xi:\Sigma \rightarrow W$. Then if $\Sigma$ is stable, $\Sigma =rW$ for some real number $r$.
\end{theorem}
The stability of the Wulff shape $W$ follows from Wulff's Theorem which states that $W$ is, in fact, the absolute minimizer for the anisotropic energy among all surfaces enclosing the same volume.

\begin{lemma}
\label{J}
Let $J=\nu\times \cdot$ denote the almost complex structure on 
$\Sigma^*$. Then there holds
\begin{equation}
\label{rep}
J\:d\xi+(d\xi)^tJ=-\Lambda \:J\:.\end{equation}
\end{lemma}
{\it Proof.}  On a neighborhood of an arbitrary point $p\in \Sigma^*$,  we can choose an orthonormal frame which diagonalizes $A:=D^2\gamma +\gamma I$.  Using this frame, we write
\begin{equation}
\label{matrices} A=\left( \begin{array}{cc}
1/\mu_1& 0 \\
0 & 1/\mu_2
\end{array} \right ) \:,\quad d\nu=\left( \begin{array}{cc}
-h_{11}& -h_{12} \\
-h_{12} & -h_{22}
\end{array} \right )  \:,\quad J=\left( \begin{array}{cc}
0& -1\\
1 & 0
\end{array} \right )\:.\end{equation}

Recalling that $d\xi=A\:d\nu$, the verification of \rf{rep} is straightforward. {\bf q. e. d. }\\[4mm]

Recall that the Jacobi operator of $\Sigma^*$ is the elliptic self-adjoint operator acting on functions given by
$$L[\psi]=\nabla \cdot A\nabla \psi +\langle d\xi, d\nu \rangle \psi\:.$$
The geometric meaning of this operator is that for a variation of the surface $\delta X= V+\psi \nu$ with $V$ tangent, the pointwise variation of $\Lambda$ is
\begin{equation}
\label{deltaL}
\delta \Lambda =L[\psi]+\nabla \Lambda \cdot V\:.
\end{equation}
\begin{lemma}
\label{delta}
On $\Sigma^*$ there holds
\begin{equation}
L[\gamma]= \Lambda^2-2K_\Sigma /K_W\:,\end{equation}
\end{lemma}
{\it Proof} We consider the variation of $\Sigma^*$ given by
\begin{equation}
\label{thevariation}
X+\epsilon \xi=X+\epsilon (D\gamma +\gamma \nu)\:.
\end{equation} 

This variation is parallel in the sense that it does not change the tangent plane to the surface. It follows that both $\gamma(\nu)$ and
$\xi=\nabla {\tilde \gamma}(\nu)$ are preserved by the variation.

Near a point $\in {\Sigma}$, let $e_1, e_2$ be a fixed orthonormal frame. By the remark above and  \rf{CH}, along the variation we can write
$$-\Lambda(\epsilon)=g^{ij}X_i\xi_j\:,$$
where $g^{ij}$ and $X$ depend on $\epsilon$ and the subscripts denote differentiation with respect to the framing vectors. It is important to note that here we are only taking the divergence along the 
surface  as opposed to using the divergence in ${\bf R}^3$. This is justified by using the reasoning of section 1.7 of  \cite{G}. Since ${\tilde \gamma}$ is homogeneous of degree one, $\xi$ is homogeneous  zero and so  $\nabla_\nu \xi(\nu)=0$ holds.

Letting `dot' denote differentiation with respect to $\epsilon$, we have, using the summation convention
$$-{\dot \Lambda} ={\dot g}^{ij}X_i\xi_j+g^{ij}{\dot X}_i\xi_j\:.$$
Setting $\epsilon =0$, we find that since the frame is orthonormal for the original metric,
\begin{eqnarray*}-{\dot \Lambda}&=&- {\dot g}_{ij}e_i\xi_j+g^{ij}{\xi}_i\xi_j\\
&=&-(\xi_i \cdot e_j+\xi_j \cdot e_i)e_i\cdot \xi_j+\langle d\xi, d\xi \rangle\:.\end{eqnarray*}
If we write $d\xi=(d\xi)_{\cal S}+(d\xi)_{\cal A}$, where the two terms are respectively the symmetric and skew symmetric parts of $d\xi$, then the previous equation reduces to
$$\dot \Lambda=\langle (d\xi)_{\cal S},(d\xi)_{\cal S} \rangle -
\langle (d\xi)_{\cal A},(d\xi)_{\cal A}\rangle\:.$$

By using the frame and matrices in \rf{matrices}, we find
$$(d\xi)_{\cal S}=\left( \begin{array}{cc}
\frac{-h_{11}}{\mu_1}&\frac{ -h_{12} }{2}(\frac{1}{\mu_1}+\frac{1}{\mu_2})\\
\frac{ -h_{12} }{2}(\frac{1}{\mu_1}+\frac{1}{\mu_2})& \frac{-h_{22}}{\mu_2}
\end{array} \right )\:,$$ 
$$ (d\xi)_{\cal A}=\left( \begin{array}{cc}
0&\frac{ -h_{12} }{2}(\frac{1}{\mu_1}-\frac{1}{\mu_2})\\
\frac{ -h_{12} }{2}(\frac{1}{\mu_2}-\frac{1}{\mu_1})& 0
\end{array} \right )  \:.$$
By \rf{deltaL} and the last three equations, we obtain
\begin{eqnarray*}
L[\gamma]&=& \langle (d\xi)_{\cal S},(d\xi)_{\cal S} \rangle -\langle (d\xi)_{\cal A},(d\xi)_{\cal A}\rangle\\
&=& (\frac{h_{11}}{\mu_1})^2+(\frac{h_{22}}{\mu_2})^2 +(\frac{h^2_{12}}{2})\bigl( [\frac{1}{\mu_1}+\frac{1}{\mu_2}]^2- [\frac{1}{\mu_1}-\frac{1}{\mu_2}]^2\bigr)\\
&=& (\frac{h_{11}}{\mu_1})^2+(\frac{h_{22}}{\mu_2})^2 +2(\frac{h^2_{12}}{\mu_1 \mu_2})\\
&=& \biggl( \frac{h_{11}}{\mu_1}+\frac{h_{22}}{\mu_2}\biggr)^2 -2\biggl( \frac{h_{11}h_{22}-h_{12}^2}{\mu_1\mu_2}\biggr)\\
&=& \Lambda^2-2K_\Sigma /K_W\:.
\end{eqnarray*}
\begin{lemma} 
\label{int}
On $\Sigma^*$, there holds
\begin{equation}
\label{div1}
\nabla\cdot J(X \times \xi)^T =2\gamma+\Lambda q \:,
\end{equation}
\begin{equation}
\label{div2}
\nabla \cdot  (d\xi+\Lambda \:I)J(X \times \xi)^T=\frac{2qK_\Sigma}{K_W}+\Lambda \gamma\:,
\end{equation}
where the superscript $T$ denotes the tangential part.
\end{lemma}
{\bf Remark} These formulas are special cases of those appearing in the work of He and Li \cite{HeLi}.\\
{\it Proof.} Writing $X=q\nu +\nabla X^2/2$, $\xi=\gamma \nu+D\gamma$, one obtains $(X\times \xi)^T= qJD\gamma -\gamma J \nabla X^2/2$, so
$J(X\times \xi)^T= \gamma \nabla (X^2/2)- qD\gamma$. 

If $h$ is any smooth function on $S^2$ with gradient $Dh$ on $S^2$ and $h_1:=h\circ \nu$ is the pull back of $h$ via the Gauss map of a smooth surface, then it is easily verified
that $\nabla h_1= d\nu(Dh)$. Here we are identifying $Dh$ with a locally defined vector field on the surface by parallel translation.

Using this, we get
$$\nabla\cdot J(X \times \xi)^T = \gamma \Delta (X^2/2)-qD\gamma +\nabla \gamma \cdot \nabla (X^2/2)-\nabla q\cdot D\gamma\:.$$
We use that $\Delta (X^2/2)=2(1+Hq)$, $\Lambda =-\nabla \cdot D\gamma +2H\gamma$ (from \rf{CH}) and $\nabla q=d\nu \nabla (X^2/2)$, together with the remarks above, to get \rf{div1}.
\begin{lemma}
\label{l2}
Under the assumptions of Theorem \rf{th1},
\begin{equation}
\label{I1}
\int_\Sigma 2\gamma +\Lambda q\:d\Sigma =0\:,\end{equation}
and 
\begin{equation}
\label{I2}\int_\Sigma \frac{2qK_\Sigma}{K_W}+\Lambda \gamma\:d\Sigma =0\:.\end{equation}
\end{lemma}
{\it Proof.}  Let $t$ and $n$ be respectively the unit tangent and the unit normal at smooth points of a $\Gamma_i$. By \rf{div1}, we have

\begin{eqnarray*} 
\int_\Sigma 2\gamma +\Lambda q\:d\Sigma&=&\sum \int_{\Sigma_j }
\nabla\cdot J(X \times \xi)^T \:d\Sigma \\
&=&\sum_j \oint_{\partial \Sigma_j}  J(X \times \xi)^T \cdot n\:ds\\
&=&-\sum_j \oint_{\partial \Sigma_j}  (X \times \xi)^T \cdot t\:ds
\end{eqnarray*}
In the last sum, each arc $\Gamma_i$ of a $\partial \Sigma_j$ is traversed twice, with opposite orientations.  So by continuity of both $X$ and $\xi$ on all of $\Sigma$, the last sum vanishes.

In a similar way, we have from \rf{div2}
\begin{equation} 
\int_\Sigma  \frac{2qK_\Sigma}{K_W}+\Lambda \gamma \:d\Sigma=\sum_j\oint_{\partial \Sigma_j}
(d\xi+\Lambda \:I)J(X \times \xi)^T \cdot n\:ds \:.
\end{equation}
Using Lemma \rf{J}, we find $J(d\xi)^t+\Lambda J=-d\xi\:J$ and so
\begin{equation*}
(d\xi+\Lambda \:I)J(X \times \xi) ^T\cdot n =-(X \times \xi) \cdot J( [d\xi]^t+\Lambda I)n =(X\times \xi)\cdot d\xi(t)\:.
\end{equation*}
Again, because $\xi$ is continuous across each $\Gamma_i$ and $\xi$ is differentiable along each $\Gamma_i$, the integrals of this form over $\Gamma_i$ considered with its two
orientations will cancel and \rf{I2} follows. {\bf q.e.d.}

\begin{prop}
\label{dVdF}
For the variation \rf{thevariation}, there holds
\begin{eqnarray*}
dV_\epsilon&=&\frac{1}{3} \biggl( q+\epsilon [\gamma -\Lambda q]+\epsilon^2 [\frac{qK_\Sigma}{K_W}-\Lambda \gamma ]+\epsilon^3 
\:\frac{\gamma K_\Sigma}{K_W}\:\biggl)\:d\Sigma\\
&=&\frac{1}{3} \biggl( q+\epsilon [\frac{-3\Lambda q}{2}+\frac{1}{2}
\nabla \cdot J(X\times \xi )^T]+\epsilon^2 [\frac{3\Lambda^2 q}{4}+\frac{1}{2} \nabla \cdot  (d\xi+\Lambda \:I)J(X\times \xi )^T -\frac{3\Lambda}{4}\nabla \cdot J(X\times \xi )^T]\\
&&+\epsilon^3 
\:\frac{\gamma K_\Sigma}{K_W}\:\biggl)\:d\Sigma\\ 
\end{eqnarray*}
and
\begin{equation}
\label{dF}
d{\cal F}_\epsilon :=\gamma\: d\Sigma_\epsilon=\gamma \biggl(1-\epsilon \Lambda +\epsilon^2 \frac{K_\Sigma}{K_W}\biggr)\:d\Sigma\:.
\end{equation}
\end{prop}
{\it Proof.} Let $e_i$ be a positively oriented  orthonormal frame at a point $p\in \Sigma^*$. We have $d\Sigma_\epsilon =||dX_\epsilon (e_1)\times dX_\epsilon(e_2)||\: d\Sigma_0=||(I+\epsilon d\xi)(e_1)\times (I+\epsilon d\xi)(e_2)||\:d\Sigma_0
=(1-\epsilon \Lambda +\epsilon^2 K_\Sigma /K_W)\:d\Sigma_0 $. Multiplying this by $\gamma$ gives \rf{dF}. Taking the product of $dX_\epsilon (e_1)\times dX_\epsilon(e_2)=(1-\epsilon \Lambda +\epsilon^2 K_\Sigma /K_W) \nu$ with $X+\epsilon \xi$ gives the first line in the equation for $dV_\epsilon$. The second line follows by Lemma \rf{int}.
{\bf q.e.d.}\\[4mm]
{\it Proof of Theorem \rf{th1}}. We now use the method of \cite{P} which was adapted from \cite{W}.
Let $X_\epsilon:=X+\epsilon \xi$ and consider the variation $s(\epsilon) (X_\epsilon)$, where the function $s(\epsilon)$ is determined so that the variation is volume preserving, i. e. 
\begin{equation}
\label{V}
V(\Sigma) \equiv V(s(\epsilon) (X_\epsilon ))=(s(\epsilon))^3
 V( X_\epsilon )\:.\end{equation}
 Write $V(X_\epsilon)=:v_0+v_1\epsilon +v_2\epsilon^2 
 +v_3\epsilon^3$ and write $s(\epsilon)=:s_0+s_1\epsilon +s_2\epsilon^2+...$ It is easily established by computing the first two derivatives of $(s(\epsilon))^3V( X_\epsilon )$ at $\epsilon=0$, that
 $$s_0=1\:,\quad s_1=\frac{-v_1}{3v_0}\:,\quad s_2=\frac{2}{9}(\:\frac{v_1}{v_0}\:)^2-\frac{1}{3}\:\frac{v_2}{v_0}\:\:.$$
 
By using the second equation of Proposition \rf{dVdF} , Lemma \rf{l2}, \rf{div1} and \rf{div2}, we find
\begin{equation}
\label{vs}
\frac{v_1}{v_0}=\frac{-3\Lambda}{2}\:, \quad \frac{v_2}{v_0}=\frac{3\Lambda^2}{4}\:.
\end{equation}
 By \rf{dF}, the energy density of $s(\epsilon)X_\epsilon$ will be
$$d{\cal F}[s(\epsilon)X_\epsilon]=(s(\epsilon))^2 {\cal F}[X_\epsilon]= \biggl(1+2s_1\epsilon +(s_1^2+2s_2)\epsilon^2+...\biggr)\gamma \biggl(1-\epsilon \Lambda +\epsilon^2 \frac{K_\Sigma}{K_W}\biggr)\:d\Sigma\:.$$
Isolating the $\epsilon^2$ coefficient, we find
$$\delta^2 {\cal F}=\partial^2_{\epsilon \epsilon}\biggl( {\cal F}[s(\epsilon)X_\epsilon] \biggr)_{\epsilon -=0}=\int_{\Sigma^*} \gamma \biggl( \frac{K_\Sigma}{K_W}-\frac{\Lambda^2}{4}\biggr)\:d\Sigma\:.$$
For the surface to be stable, the integrand must be non negative. However 
$$\frac{\Lambda^2}{4}-\frac{K_\Sigma }{K_W}= \frac{1}{4}\biggl(\frac{h_{11}}{\mu_1}+\frac{h_{22}}{\mu_2} \biggl)^2-\frac{h_{11}h_{22}-h_{12}^2}{\mu_1\mu_2} =\frac{1}{4}\biggl(\frac{h_{11}}{\mu_1}-\frac{h_{22}}{\mu_2}\biggr)^2 +\frac{h_{12}^2}{\mu_1\mu_2}\:.$$
Thus, stability implies that  $(1/4)\Lambda^2-K_\Sigma/K_w \equiv 0$ holds.  Multiplying this equation by $q$ and integrating using \rf{I2} leads to
$$-\frac{\Lambda}{2} {\cal  F}[\Sigma]=\frac{3\Lambda^2}{4}V[\Sigma]\:.$$
While multiplying by $\gamma$ and integrating, using \rf{I1}, then gives
\begin{equation}
\label{in}
{\rm deg}(\xi){\cal F}[W]=\int_\Sigma \frac{\gamma K_\Sigma}{K_W}\:d\Sigma =\frac{\Lambda^2}{4}{\cal F}[\Sigma] =\frac{({\cal F}[\Sigma])^3}{9(V[\Sigma])^2}\:,\end{equation}
by the previous equation. 

We claim that  deg $\xi\le 1$ holds.  Since $\Sigma$ and $W$ are piecewise smooth, they can be continously deformed to a smooth surface ${\tilde \Sigma}$ and a smooth convex surface ${\tilde
W}$ within an $\epsilon$ neighborhood of each surface in the $C^0$ topology.  The Cahn-Hoffman map ${\tilde \xi}:{\tilde \Sigma}\rightarrow {\tilde W}$ has the same degree as
$\xi$.  The map  ${\tilde \xi}$ factors ${\tilde \xi}=\chi \circ {\tilde \nu}$ where ${\tilde \nu}$ is the Gauss map of ${\tilde \Sigma}$ and $\chi:S^2 \rightarrow {\tilde \Sigma}$ is the inverse
of the Gauss map of ${\tilde W}$. Thus $\chi$ is a bijection, since ${\tilde W}$ is convex, and so deg$\: {\tilde \xi}=$deg$\: {\tilde \nu}\le 1$ by the Gauss Bonnet Theorem.

To finish the proof, note that \rf{in} now reads 
$$\frac{({\cal F}[W])^3}{(V[W])^2}= 9{\cal F}[W]\ge \frac{({\cal F}[\Sigma])^3}{(V[\Sigma])^2}\:,$$
A strict inequality in the previous line would violate the  the isoperimetric inequality for the functional ${\cal F}$, \cite{BM}.  Also, any surface which realizes equality in this isoperimetric inequality   is a rescaling of $W$ up to a set of zero measure \cite{BM}. However both $\Sigma$ and $W$ are of class $C^0$ so, in fact, $\Sigma=rW$ for some $r>0$. {\bf q.e.d}\\[4mm]
{\bf Remark.} One can use this result to treat the case where $W$ can have point singularities of the type that would arise if $W$ were axially symmetric and was generated by rotating an arc of a circle making an acute angle with the rotation axis. If $\Sigma$ is an equilibrium surface having a non smooth points where the Cahn-Hoffman field extends
continuously, then `edges' made up of smooth curves on $\Sigma$ can be introduced and the previous results still apply. \\[4mm]
\begin{figure}
  \hfill
  \begin{minipage}[n]{.45\textwidth}
    \begin{center}
       \framebox{\includegraphics[width=60mm,height=60mm,angle=0]{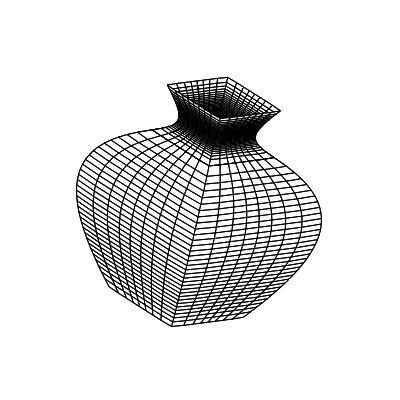}}
    \end{center}
     \caption{An unduloid for the Wulff shape in Figure 1}
 \end{minipage}
\hfill


    \end{figure}
    The construction of anisotropic Delaunay surfaces found in \cite{KP1} and \cite{KP2} can be extended to the case where the Wulff shape has edges. We begin with the assumption that $W$ is of `product form', which means that $W$ can be parameterized $\chi(\sigma, \tau)=(u(\sigma)(\alpha(\tau), \beta(\tau), v(\sigma))$, where $(u,v)$ and $(\alpha, \beta)$ are convex plane curves. We then seek a surface $\Sigma$ which can be parameterized $X(s,\tau)=(x(s)(\alpha(\tau), \beta(\tau),z(s))$. This means that the cross sections of $\Sigma$ are all homothetic to the $(\alpha, \beta)$ curve. In the case where $W$ is smooth, it was shown that the condition that $\Sigma$ have $\Lambda \equiv$ constant, is independent of the choice of the cross sectional curve $(\alpha, \beta)$. Thus, to find the $(x,z)$ curve, one could just assume that the $(\alpha, \beta)$ curve is a circle and use the results of \cite{KP1} where all such curves were found and the resulting surfaces were completely classified.
    
    If the $(\alpha, \beta)$ curve is continuous but only piecewise smooth,  this construction is equally valid at all points $(s,\tau)$ which the $(\alpha, \beta)$ curve is smooth at $\tau$. At smooth points of $\Sigma$, the Cahn-Hoffman map $\xi$ is given in the local coordinates by $(s,\tau)\mapsto (\sigma, \tau)$ by the conditions that the tangent planes to $\Sigma$ and $W$ agree at the corresponding points.   By the remarks above, $u$ and $v$ are globally smooth functions of $s$, since this is the case when the cross sections are changed to circles. Therefore, $\xi$ is globally continuous for these surfaces. Examples are shown in Figures 1 and 2.

\begin{flushleft}
Bennett P{\footnotesize ALMER}\\
Department of Mathematics\\
Idaho State University\\
Pocatello, ID 83209\\
U.S.A.\\
E-mail: palmbenn@isu.edu
\end{flushleft}

\end{document}